\documentclass[11pt]{article}   

\newtheorem{theorem}{Theorem}[section]
\newtheorem{lemma}[theorem]{Lemma}
\newtheorem{definition}[theorem]{Definition}
\newtheorem{proposition}[theorem]{Proposition}
\newtheorem{corollary}[theorem]{Corollary}

\newtheorem{remark}[theorem]{Remark}
\newtheorem{proof}[theorem]{Proof}
\newtheorem{example}[theorem]{Example}

\usepackage{amssymb,amsmath,mathrsfs}

\usepackage{pstricks,pst-node,pst-plot}

\usepackage[english]{babel}

\def\proof{{\noindent \bf proof :} }
\def\endproof{{\hfill $\square$} \vspace{2 mm}

}

\def\og{\leavevmode\raise.3ex\hbox{$\scriptscriptstyle\langle\!\langle$~}}
\def\fg{\leavevmode\raise.3ex\hbox{~$\!\scriptscriptstyle\,\rangle\!\rangle$}}

\def\ntilde{$\tilde{\rm n}$}

\def\ccC{\mathscr{C}}

\def\ccG{\mathscr{G}}
\def\ccK{\mathscr{K}}
\def\ccP{\mathscr{P}}

\def\ms{\medskip}
\def\ss{\smallskip}

\def\s{\sigma}
\def\udl{\underline}
\def\red{\searrow \!\!\stackrel{s}{}}
\def\exp{\stackrel{s}{} \!\!\nearrow }
\def\redws{\searrow \!\!\stackrel{ws}{}}

\def\redi{\searrow \!\!\stackrel{I}{}}
\def\expi{\stackrel{I}{} \!\!\nearrow }

\usepackage{color}

\begin{document}

\title{Simplicial simple-homotopy of flag complexes in terms of graphs}
\author{
Romain Boulet, 
Etienne Fieux
and 
Bertrand Jouve  }
\date{ Institut de Math\'ematiques de Toulouse \\ 
Unit\'e Mixte de Recherche \\ CNRS - UMR 5219
}
\maketitle

\abstract{A flag complex can be defined as 
a simplicial complex whose simplices correspond to complete subgraphs
of its 1-skeleton taken as a graph. In this article, by introducing 
the notion of s-dismantlability,
 we shall define the s-homotopy type of a graph and show in particular  
that two finite graphs have the same s-homotopy type 
if, and only if, the two flag complexes determined by these graphs have 
the same simplicial simple-homotopy type
(Theorem \ref{theoremGK}, part 1).
This result is closely related  to similar results established by 
Barmak and Minian (\cite{barmin}) in the framework of posets
and we give the relation between the two approaches
(theorems \ref{theoremGP} and \ref{theoremPK}). 
We conclude with a question about the relation between the s-homotopy
and the graph homotopy defined in \cite{cyy}. }
\vspace{ 2 mm}


\vspace{3 mm}

\noindent {\large \bf Introduction}
\vspace{2 mm}

Flag complexes are (abtract) simplicial complexes whose every minimal 
non simplex has two elements 
(\cite{stanley},\cite{renteln02},\cite{frohmader}); 
this means that a flag complex is completely determined 
by its 1-skeleton (all necessary definitions are recalled below). 
They constitute an important subset of the set of simplicial complexes;
in particular, the barycentric subdivision of any simplicial complex 
is a flag complex and we know that a simplicial complex
and its barycentric subdivision have the same simple-homotopy type.

Flag complexes arise naturally from the graph point of  view and are also 
sometimes called {\sl clique complexes} (\cite{cyy}).
Indeed, the 1-skeleton of a simplicial complex can be considered as a graph
and it is easy to see that  a simplicial complex  $K$ is a flag complex
if, and only if,
we can write $K=\Delta_{\ccG}(G)$ for some graph $G$ 
where, by definition, $\Delta_{\ccG}(G)$ 
is the simplicial complex whose simplices are given by the 
complete subgraphs of $G$ (and this is sometimes
taken as the definition of flag complexes, as in \cite{kozlov_CAT}).

In this paper, we are interested in the notion of simplicial simple-homotopy 
 for flag complexes.
We note that the determination of the simplicial simple-homotopy type 
 is actually important not only for simplicial complexes 
but also for graphs because simplicial complexes arise in various constructions  
in graph theory. 
For example, this notion appears in the study of the clique graph 
(\cite{prisner92}, \cite{lapivf08}),
in results about the polyhedral complex ${\textsf{Hom}}(G,H)$ 
introduced by Lovasz (\cite{babsonkozlov},\cite{kozlov2}) or in relation to 
evasiveness (\cite{kss}).

Simplicial simple-homotopy is defined by formal deformations
themselves defined by the notion of elementary collapses consisting 
in the deletion of certain
pairs of simplices (see \S 2). As the set of simplices of a flag complex 
is determined by its 1-skeleton seen as a graph,
the aim of this paper is to relate formal deformations on flag complexes 
to certain operations on graphs.
The key notion will be the one of {\sl s-dismantlability}: the deletion or the addition 
of  {\sl s-dismantlable} vertices in a graph will play the role of elementary reductions
or expansions 
(\cite{cohen}) in a simplicial complex.
More precisely, a vertex $g$ of a graph $G$ will be called s-dismantlable 
if its open neighborhood
is a dismantlable graph and  the deletion of an s-dimantlable vertex $g$ in $G$
is equivalent to the deletion of all simplices which contain $g$ in 
$\Delta_{\ccG}(G)$.

In Section \ref{sec_s-dism}, we introduce  the notion of s-dismantlability which allows us 
to define an equivalence relation
for graphs; the equivalence class $[G]_s$ of a graph $G$ for this equivalence relation
will be called the {\sl s-homotopy type of $G$} and we  give some 
properties related to these notions.

In Section \ref{sec_complexe},   we study the correspondence  between s-dismantlability in $\ccG$ 
(the set of finite undirected graphs, without multiple edges) 
and simplicial simple-homotopy
 in  $\ccK$ (the set of  finite simplicial complexes);
we prove that two finite graphs $G$ and $H$
have the same s-homotopy type if, and only if, $\Delta_{\ccG}(G)$ and $\Delta_{\ccG}(H)$ have the same
simple-homotopy type.
  Reciprocally, for a simplicial complex $K$, $\Gamma(K)$ is the graph whose vertices
are the simplices of $K$ and the edges are given by the inclusions; 
we show that two finite simplicial
complexes $K$ and $L$ have the same simple-homotopy type 
if, and only if, $\Gamma(K)$
and $\Gamma(L)$ have the same s-homotopy type.
We have to mention that these results need the introduction of  
barycentric subdivision in $\ccG$, defined, for a graph $G$, as 
the 1-skeleton of  the usual barycentric subdivision (in $\ccK$) of $\Delta_{\ccG}(G)$.

In Section \ref{sec_poset}, we consider the  important class of flag complexes   
which results from posets: if $P$ is a poset,
 $\Delta_{\ccP}(P)$ is the flag complex whose simplices are given by chains of $P$.
In (\cite{barmin}), Barmak and Minian define a notion of \emph{simple equivalence} 
in posets
and show that there is a one-to-one correspondence between simple-homotopy 
types of finite simplicial
complexes and simple equivalence classes of finite posets. 
As we have $\Delta_{\ccP}(P)=\Delta_{\ccG}(Comp(P))$
where $Comp(P)$ is the comparability graph of $P$,
there is a close relation between this approach and our approach from graphs. 
 We show that there is indeed
a one-to-one correspondence between s-homotopy type in  $\ccG$ 
and simple equivalence classes in the set $\ccP$ of finite posets.
Finally, we consider a triangle between finite graphs, posets and simplicial complexes
recapitulating the close relations between s-homotopy type (in $\ccG$), 
simple type (in $\ccP$)
and simple homotopy type (in $\ccK$).

In Section \ref{sec_weak} we describe a weaker version of s-dimantlability 
on graphs which provides a closer connection with simplicial collapse 
for flag complexes (Proposition \ref{ws_prop_G_to_K}).

Then we conclude in Section \ref{sec_Chen} with a question concerning 
the relation between s-homotopy and the graph homotopy defined in \cite{cyy}

Some results  have been set out in \cite{bfj}.

\medskip

\noindent {\bf Definitions, notations} 
\medskip

Let $\mathscr{G}$ be the set of finite undirected graphs, without multiple edges. 
If $G\in \mathscr{G}$, we have  $G=(V(G),E(G))$ with 
$E(G) \subseteq \{ \:\{g,g'\}\:,~g,g'\in V(G)\}$.
For brevity, we write $xy \in G$ or  $x\sim y$  for $\{x,y\} \in E(G)$ 
and $x \in G$ for $x \in V(G)$.
The closed neighborhood of  $g$ is  $N_G[g]:=\{ h \in G\:,\:g \sim h\} \cup \{g\}$
and $N_G(g) :=N_G[g]\setminus\{g\}$ is its open neighborhood.
When no confusion is possible, a subset $S$ of $V(G)$ will also denote
the subgraph of $G$ induced by  $S$. We denote by $G \setminus S$ the graph 
obtained from $G$ by 
deleting $S$ and all the edges adjacent to a vertex of $S$. In particular, 
we use the notations $G \setminus x$ 
and $G \setminus xy$ to indicate the deletion of a vertex $x$ or an edge $xy$.
The notation  $c=[g_1 ,\ldots ,g_k]$ means that the the subset 
$\{g_1,\ldots,g_k\}$ of $V(G)$
induces a complete subgraph of $G$. 
We say that a graph $G$ is a {\sl cone} on a vertex $g \in G$ if $N_G[g]=G$.
The notation $pt$ will denote a graph reduced to a single vertex (looped or not looped).

Let  $\ccK$ be the set of (abstract) finite simplicial complexes. 
A simplicial complex $K$ is a family of subsets of a finite set $V(K)$
(the set of vertices of $K$) stable with respect to  deletion of  elements
(if $\s \in  K$ and $x \in \s$, then $\s \setminus \{x\} \in K$).
An element $\{x_0,x_1,\ldots,x_{k}\}$ of $K$ is called $k$-simplex and
will be denoted by $< x_0,x_1,\ldots,x_{k}>$;
the   $n$-skeleton of $K$ is the set $K_n$ formed by all 
$k$-simplices of $K$ with $k\leq n$.
 The 0-skeleton is identified with the set $V(K)$ of vertices of  $K$. 
A face of $\sigma =< x_0,x_1,\ldots,x_{k}>$ 
is any simplex included in $\sigma$.

\section{s-dismantlability and s-homotopy type in $\ccG$}\label{sec_s-dism}

\subsection{Definitions}

Let $G \in \mathscr{G}$. We recall (\cite{prisner92},\cite{bcf94},\cite{ginsburg}) 
that a vertex  $g$ of $G$ is called
{\sl dismantlable} if  there is another vertex  $g'$ of $G$ which dominates $g$
 (\textit{i.e.}, $g \neq g'$ and $N_G[g] \subseteq N_G[g']$); 
we note that this implies $g\sim g'$. 
A graph  $G$ is called {\sl dismantlable} if it is reduced to a single vertex or if we can write
$V(G)=\{g_1,g_2,\ldots,g_n\}$ with $g_i$ dismantlable in the subgraph induced 
by $\{g_1,\ldots,g_i\}$, for $2\leq i \leq n$. In particular, a dismantlable graph is necessarily
non empty.
\ss

\begin{definition} A vertex $g$ of a graph  $G$ is called s-dismantlable  in $G$ if $N_G(g)$ is dismantlable. 
\end{definition}

Let $H$ a subgraph of a graph $G$. We shall say that \textit{$G$ is dismantlable on $H$} if we can go from $G$ to $H$ by 
successive deletions of dismantlable vertices.

\begin{example} \label{minimal_s_dismantlable}
No vertex of the graph $G_1$ of Fig. \ref{8_vertex_s_collapsible} is 
dismantlable but there are four s-dismantlable vertices, 
as  the vertex $a$ in the picture $(N_{G_1}(a)$ is a dismantlable path).

\begin{figure}[h]
\psset{xunit=0.5cm,yunit=0.5cm}
\pspicture(-10,-2.5)(2,3.5)
\rput(-1,-1){$\bullet$}
\rput(1,1){$\bullet$}
\rput(-1,1){$\bullet$}
\rput(1,-1){$\bullet$}
\rput(-2,-2){$\bullet$}
\rput(-2,2){$\bullet$}
\rput(2,-2){$\bullet$}
\rput(2,2){$\bullet$}
\psline(1,1)(-1,1)(-1,-1)(1,-1)(1,1)(-1,-1)(-2,2)(-1,1)(2,2)(1,1)(2,-2)(1,-1)(-2,-2)(-2,2)(2,2)(2,-2)(-2,-2)(-1,-1)
\rput(-4,0){$G_1$}
\endpspicture
\psset{linecolor=lightgray}
\pspicture(-9,-2.5)(2,3.5)
\rput(-1,-1){$\bullet$}
\rput(1,1){$\bullet$}
\rput(-1,1){$\bullet$}
\rput(1,-1){$\bullet$}
\rput(-2,-2){$\bullet$}
\rput(-2,2){$\bullet$}
\rput(2,-2){$\bullet$}
\rput(2,2){$\bullet$}
\psline(1,1)(-1,1)(-1,-1)(1,-1)(1,1)(-1,-1)(-2,2)(-1,1)(2,2)(1,1)(2,-2)(1,-1)(-2,-2)(-2,2)(2,2)(2,-2)(-2,-2)(-1,-1)
\uput[ur](2,2){$a$}
\psline[linecolor=black](-2,2)(-1,1)(1,1)(2,-2)
\rput(5,0){$N_{G_1}(a)$}
\endpspicture
\caption{The vertex $a$ is s-dismantlable (and not dismantlable)}
\label{8_vertex_s_collapsible}
\end{figure}

\end{example}

\noindent Clearly, a dismantlable vertex  is s-dismantable 
(because the open neighborhood
of a dismant\-la\-ble vertex is a cone which is a dismantlable graph). 
When $g \in G$ is s-dismantlable, we shall write $G \red G\setminus \{g\}$ (elementary reduction); equivalently,
 $G \exp G \cup \{x\}$ (elementary expansion) indicates the addition to  $G$ of a vertex $x$
such that $N_{G \cup \{x\}}(x)$ is dismantlable.
By analogy with the usual situation in $\ccK$,   $G \red H$ (resp. $G \exp H$) indicates that we can go
from  $G$ to $H$,
by deleting (resp. adding) successively  s-dismantlable vertices. 

\begin{definition} \label{def_s_homotopy_type}
Two graphs $G$ and $H$ have the same s-homotopy type
if there is a sequence  $G=J_1,\ldots,J_k=H$ in  $\mathscr{G}$ such that
$G =J_1\stackrel{s}{\to} J_2 \stackrel{s}{\to}  \ldots \stackrel{s}{\to} J_{k-1} 
\stackrel{s}{\to} J_k=H$ 
where each arrow $\stackrel{s}{\to} $ represents the suppression or the addition 
of an s-dismantlable vertex.
\end{definition}

\noindent  This defines an equivalence relation in $\mathscr{G}$ and we shall 
denote by $[G]_s$
the equivalence class representing the s-homotopy type of a graph $G$.
A graph $G$ will be called s-dismantlable if $[G]_s=[pt]_s$.

\subsection{Properties}

\begin{lemma} \label{lemma_G_W_H}
Let $G,H \in \mathscr{G}$ such that $[G]_s=[H]_s$; then, there is a graph $W$
such that $G \exp W \red H$.
\end{lemma}

\proof Let us suppose that $G$ and $H$ are two graphs such that $[G]_s=[H]_s$. It is sufficient
to prove that an elementary reduction and an elementary expansion may be switched in the sequence
of elementary operations from $G$ to $H$. So, let us suppose that in a graph
$G'$, we have an elementary reduction followed by an elementary expansion:
 $$ (1)~~~~~ G' \red G' \setminus g_1 \exp (G' \setminus g_1) \cup \{g_2\} $$
This means that the graphs $N_{G'}(g_1)$ and $N_{G'\setminus g_1}(g_2)$ 
are dismantlable and, in
particular, $g_2 \not\sim g_1$. So we can adjoin $g_2$ to $G'$ 
by putting $N_{G'}(g_2)=N_{G'\setminus g_1}(g_2)$;
of course, $N_{G'\cup\{g_2\} }(g_1)=N_{G'}(g_1)$ and the sequence $(1)$ can be 
alternatively written
 $$(2)~~~~~ G' \exp G' \cup \{g_2\} \red (G' \cup \{g_2\}) \setminus g_1$$ 
with isomorphic resultant  graphs
$ (G' \setminus g_1) \cup \{g_2\}$ and $(G' \cup \{g_2\}) \setminus g_1$ and 
this proves that all the reductions can be 
pushed at the end of the sequence of elementary operations  from $G$ to $H$.
\endproof

We recall that the suspension  $SG$ of $G$ is the graph whose vertex 
set is $V(G)\cup \{x,y\}$ where $x$ and $y$ are 
two distinct vertices which are not in $V(G)$ 
and whose edge set is $E(G)\cup \{\: xg,g \in V(G)\} \cup \{\:yg,g \in V(G)\} $; 
in the sequel, 
$SG$ will be also denoted by $G\cup \{x,y\}$.
We shall need the following result :

\begin{proposition} \label{G_SG_dismantlable}
A graph $G$ is dismantlable if, and only if, its suspension $SG$ is dismantlable.
\end{proposition}

\proof
Let $SG=G\cup \{x,y\}$. If $G$ is dismantlable, let 
$V(G)=\{g_1,g_2,\ldots,g_n\}$ with $g_i$ dismantlable in the subgraph induced 
by $\{g_1,\ldots,g_i\}$, for $2\leq i \leq n$. 
Then we can write $V(SG)=\{g_1,x,y,g_2,\ldots,g_n\}$ with $g_i$ 
dismantlable in the subgraph induced  by $\{g_1,x,y,\ldots,g_i\}$, 
for $2\leq i \leq n$ and the subgraph of $SG$
induced by $\{g_1,x,y\}$ is of course dismantlable (it is a path).
Let us now suppose that  $SG$ is dismantlable and let $g$ be a dismantlable 
vertex in $SG$. 
If $g=x$ or $g=y$, this means that $G$ is a cone (because $N_{SG}(x)
=N_{SG}(y)=G$ and a vertex $g'$ which dominates $g$ verifies $N_G[g']=G$). 
If $g\neq x$ and $g \neq y$, $g$ is also dismantlable in $G$ 
(because we have $\{x,y\} \in N_{SG}(g)$ and a vertex which dominates 
$g$ in $SG$ is necessarily different from $x$ and $y$
and, consequently, dominates $g$ in $G$).
From this observation, it follows that when we delete a dismantlable 
vertex $g$  in $SG$, either $g\in G$ and $g$ is 
dismantlable in $G$, either $g\in \{x,y\}$ and this implies that $G$ is 
dismantlable (because it is a cone). By iteration
of this procedure, we get that $G$ is a dismantlable graph.
\endproof

\begin{lemma} \label{edge_s_deletable} (deletion of an edge) 
If $g$ and $g'$ are two distinct vertices of a graph  $G$ such that $g \sim g'$ and
 $N_G(g)\cap N_G(g')$ is nonempty and dismantlable, then $[G]_s=[G\setminus  gg']_s$. 
In other words, we can s-delete the edge $gg'$. 
\end{lemma}

\proof
 We add to $G$ a vertex $x$ with edges $xz$ for every $z$ in $N_G[g]\setminus \{g'\}$
(it is an elementary expansion because $N_G[g]\setminus\{g'\}$ is  a cone on $g$)
and we  write $G \cup x$ for the resulting graph.  
Let us  verify that $g$ is s-dismantlable
in  $G \cup x$. We have $N_{ G \cup x}(g)=N_G(g) \cup \{x\}$. 
If $N_G(g)\subseteq N_G[g']$,
we can write $N_{ G \cup x}(g)=\left(N_G(g)\cap N_G(g')\right) \cup x \cup g'$. If
$N_G(g)\not\subseteq N_G[g']$, every $y$ in $N_G(g)$ which is not in $N_G[g']$
is dominated by $x$ in $N_{ G \cup x}(g)$; so, $N_{ G \cup x}(g)$ is dismantlable on its 
subgraph  induced by the set of vertices
$\left(N_G(g)\cap N_G(g')\right) \cup x \cup g'$.
But $\left(N_G(g)\cap N_G(g')\right) \cup x\cup g'$ 
  is the suspension of $N_G(g)\cap N_G(g')$ and, by 
Proposition \ref{G_SG_dismantlable},
it is a dismantlable graph because 
$N_G(g)\cap N_G(g')$ is dismantlable.
Thus,  $g$ is s-dismantlable in $G\cup x$ and we can reduce
$G\cup x$ on $(G\cup x)\setminus g$ which 
is clearly isomorphic to $G\setminus  gg'$.
\endproof

\begin{proposition} \label{s_dismantlable_neighbor} Let $G \in \ccG$ and $g \in V(G)$.
 If $N_G(g)$ is s-dismantlable, then $[G]_s=[G\setminus g]_s$.
\end{proposition}

\proof
Let $g \in V(G)$ such that the graph $N_G(g)$ is s-dismantlable. 
 Let us suppose first that $N_G(g) \red pt$ and let $N_G(g)=\{y_1,y_2,\ldots,y_k\}$ 
with $y_i$ s-dismantlable in the subgraph of $N_G(g)$ induced  by $\{y_1,\ldots,y_i\}$,
for $2 \leq i \leq k$. 
Then, $N_G(y_k) \cap N_G(g)$ is dismantlable (it is the open neighborhood of $y_k$ in 
$N_G(g)$) and, by  Lemma \ref{edge_s_deletable}, 
$G$ is s-dismantlable on the graph $H_k=G\setminus y_kg$. 
Now, $N_{H_k}(g)=\{y_1,y_2,\ldots,y_{k-1}\}$ and  
$N_{H_k}(y_{k-1}) \cap N_{H_k}(g)=N_{G}(y_{k-1}) \cap N_{H_k}(g)$ is dismantlable 
in $H_k$
and, by the same argument,
$H_k$ is s-dismantlable on $H_{k-1}=G\setminus \{y_kg,y_{k-1}g\}$.
The iteration of this procedure shows that $G$ is s-dismantlable on 
$H_{2}=G\setminus \{ y_kg,y_{k-1}g,\ldots,y_3g,y_2g\}$.
Of course, $g$ is s-dismantlable (in fact, dismantlable) in $H_2$ and this proves that 
$G$ is s-dismantlable on $G\setminus g$.
If we don't have $N_G(g) \red pt$ , we know by Lemma \ref{lemma_G_W_H} 
that there is a graph $W$ such that
$N_G(g) \exp W \red pt$. Let $\{z_1,\ldots,z_m\}$ the set of vertices of $W$ 
which are not in 
$N_G(g)$. We define $H=G \cup \{z_1,\ldots,z_m\}$ 
with $N_H(z_i)=N_{W}(z_i) \cup \{g\}$.
It is clear that $G \exp H$ and that $N_H(g)=W \red pt$. 
Thus, by the previous discussion, we know
that $H$ is s-dismantlable on  $H \setminus g$ 
but $H \setminus g \red G \setminus g$
(because of $W \red N_G(g)$) and we can conclude that $[G]_s =[ G\setminus g]_s$
(by $[G]_s=[ H]_s =[  H \setminus g]_s =[  G \setminus g]_s$).
\endproof

This result implies that the procedure of s-dismantlability can be done more rapidly
by deleting a vertex whose open neighborhood is s-dismantlable.
For example, we  immediately get the following result:

\begin{corollary} Let $G \in \ccG$. If $G$ is s-dismantlable, then $SG$  is also s-dismantlable.
\end{corollary}

\proof
As $N_{SG}(x)=N_{SG}(y)=G$ and $G$ is s-dismantlable, we obtain by 
the previous proposition that 
 $[SG]_s = [SG\setminus \{x,y\}]_s$, \textit{i.e.} $[SG]_s =[G]_s$.
\endproof

Now, by analogy with the notion of collapsibility in 
simplicial complexes (see below), we introduce the following definition:

\begin{definition}
A graph $G$ is called s-collapsible if $G\red pt$
 (\textit{i.e.} we can write 
$G=\{g_1,g_2,\ldots,g_n\}$ with $g_i$ s-dismantlable in the subgraph induced by
 $\{g_1,\ldots,g_i\}$, for $2\leq i \leq n$).
\end{definition}

\begin{remark}
Of course, a dismantlable graph is s-collapsible but the inclusion 
(of the family of dismantlable graphs in
the family of s-collapsible graphs) is strict;  for instance, the graph $G_1$ of example 
\ref{minimal_s_dismantlable} is not dismantlable 
(because there is no dismantlable vertex) 
and s-collapsible: $G_1\red G_1\setminus a$ because $a$ is 
s-dismantlable and it is easy to see that the graph
$G_1\setminus a$ is s-collapsible (even more: $G_1\setminus a$ is dismantlable).
Furthermore, it can be proved that this graph $G_1$ is minimal 
(in terms of number of vertices) in the family of  s-collapsible 
graphs and without any dismantlable vertex. 
\end{remark}

\begin{proposition} Let $G \in \ccG$. The graph $G$ is s-collapsible if, 
and only if,  $SG$  is  s-collapsible.
\end{proposition}

\proof
Let $SG=G\cup\{x,y\}$. First, we observe that every  s-dismantlable vertex $g$ in $G$
 is also s-dismantlable in $SG$ because $N_{SG}(g)=N_G(g)\cup \{x,y\}=SN_G(g)$ 
and the dismantlability 
of $SN_G(g)$ is a consequence of the dismantlability of $N_G(g)$ 
(Proposition \ref{G_SG_dismantlable}).
Suppose now that   $G \red pt$ with
$G=\{g_1,g_2,\ldots,g_n\}$ and $g_i$ s-dismantlable in the subgraph induced by
 $\{g_1,\ldots,g_i\}$, for $2\leq i \leq n$. Then, by the previous observation,
 $g_i$ is s-dismantlable in the subgraph of $SG$ induced by
 $\{g_1,x,y, g_2\ldots,g_i\}$, for $2\leq i \leq n$. 
So, we have \linebreak $SG \red \{g_1,x,y\} \red \{g_1\}$.
Reciprocally, if $SG$ is s-collapsible, we shall prove that $G$ 
is also s-collapsible by induction on the number of vertices
of $G$. If $\vert V(G)\vert=1$, there is nothing to prove 
($G$ and $SG$ are s-collapsible). Let us suppose that
the s-collapsibility of $SG$ implies the s-collapsibility of $G$ 
if $\vert V(G)\vert=n$ for $n\geq 1$ and
let $G$ with $\vert V(G)\vert=n+1$ such that $SG$ is s-collapsible. 
Let $g\in V(SG)$ an s-dismantlable vertex.
If $g\in\{x,y\}$, this means that $G$ is dismantlable 
(because $N_{SG}(x)=N_{SG}(y)=G$) and thus, s-collapsible.
So, we can assume that $g\neq x$ and $g\neq y$. 
We have $N_{SG}(g)=SN_G(g)$; so, by Proposition \ref{G_SG_dismantlable},
the dismantlability of $N_{SG}(g)$ implies the dismantlability of $N_G(g)$ 
and we have $G\red G\setminus g$.
Now, $SG\setminus g=S(G\setminus g)$ and we can apply the induction 
hypothesis to conclude that
$G\setminus g$ is s-collapsible (and the same conclusion for $G$).
\endproof

\section{Relation with simple homotopy in $\ccK$}\label{sec_complexe}

\subsection{From $\ccG$ to $\ccK$}

Let $K \in \ccK$; let us recall (\cite{cohen}) that an elementary  simplicial 
reduction  (or {\sl collapse})  in $K$ 
is the suppression of a pair of simplices
 $(\sigma,\tau)$  of $K$ such that  $\tau$ is a proper maximal face of $\sigma$
and $\tau$ is not the face of another simplex (one says that $\tau$ 
is a  free face of $K$). 

\begin{figure}[h]
\psset{unit=0.5 cm}
\pspicture(-12,-4)(3,3.5)
\pspolygon[fillstyle=solid,fillcolor=lightgray](0,0)(1,2)(2,1)(0,0)
\psline[linewidth=1.5pt](1,2)(2,1)
\pscurve (1,2)(-1,2.5)(-2,1.5)(-3,0.2)(-2,-1)(-1,-2.2)(0,-2)(1,-1.6)(1.6,-1.2)(2.2,0)(2,1)
\rput(1.7,1.7){$\tau$}
\rput(1,1){$\sigma$}
\rput(0,-3.5){\large $K$}
\psline{->}(3,0.5)(5,0.5)
\endpspicture
\pspicture(-6,-4)(3,3.5)
\psline(1,2)(0,0)(2,1)
\pscurve (1,2)(-1,2.5)(-2,1.5)(-3,0.2)(-2,-1)(-1,-2.2)(0,-2)(1,-1.6)(1.6,-1.2)(2.2,0)(2,1)
\rput(0,-3.5){\large $K\setminus \{\sigma,\tau\}$}
\endpspicture
\caption{An elementary collapse}
\end{figure}

\noindent This is denoted by
$K \red \left( K\setminus \{\sigma,\tau\}\right)$ (and called {\sl elementary collapse}
or {\sl elementary reduction}) or 
$\left( K\setminus \{\sigma,\tau\} \right)\exp K$ ({\sl elementary anticollapse}
or {\sl elementary expansion}).
More generally, a simplicial collapse $K\red L$ (resp. anticollapse $K\exp L$) 
is a succession of elementary simplicial collapses 
(resp. elementary simplicial anticollapses) 
which tranform $K$ into $L$. Collapses or anticollapses are 
called {\sl formal deformations}.
A simplicial complex $K$ is called {\sl collapsible} if $K \red pt$.
Two simplicial complexes $K$ and $L$ have the same {\sl simple-homotopy type}
if there is  (in  $\ccK$)  a finite sequence
$K =M_1\stackrel{s}{\to} M_2 \stackrel{s}{\to}  \ldots \stackrel{s}{\to} M_{k-1} 
\stackrel{s}{\to} M_k=L$ where each arrow $\stackrel{s}{\to} $ represents 
a simplicial collapse 
or a simplicial anticollapse.  
We shall denote by $[K]_s$ the simple-homotopy type of $K$.
\ss

\noindent Let us recall that for a simplex $\sigma \in K$, 
${\rm link}_K(\sigma):=\{\tau \in K \:, 
\sigma \cap \tau=\emptyset  ~{\rm and}~\sigma \cup \tau \in K\}$ and
${\rm star}^o_K(\sigma)$ is the set of simplices of $K$ containing $\sigma$.
Let us recall :
\ms

\begin{lemma}[{\cite[Lemma 2.7]{welker}}]  \label{lemma_welker}
Let $\sigma$ be a simplex of  $K\in \ccK$. If ${\rm link}_{K}(\sigma)$ is collapsible,
then $K\red K \setminus {\rm star}^o_K(\sigma)$.
\end{lemma}

We recall that the application $\Delta_{\mathscr{G}} : \mathscr{G} \to \ccK$ 
is defined in the following way :
if $G \in \mathscr{G}$, $\Delta_{\mathscr{G}}(G)$  is the simplicial complex 
whose simplices
are the complete subgraphs of $G$ (so we have $V(\Delta_{\mathscr{G}}(G))=V(G)$).
We note that if a graph $G$ does not contain any triangle (\textit{i.e.} a complete 
subgraph with three vertices), we can identify
$G$ and $\Delta_{\mathscr{G}}(G)$ (we consider $G$ either as a graph, 
or as a simplicial complex of dimension 1);
it is the case of $C_4$ in examples of Fig. \ref{Illustration_Delta_G}.

\begin{figure}[h]
  \begin{center}
\psset{unit=0.8 cm}
\pspicture(0,0)(2.5,3)
\pspolygon(0,0)(2,0)(1,2)(0,0)
\rput(0,0){$\bullet$}
\rput(2,0){$\bullet$}
\rput(1,2){$\bullet$}
\uput[dl](0,0){ $a$}
\rput(2.1,-0.1){ $b$}
\uput[u](1,2){ $c$}
\rput(1,-1.1){\large $K_3$}
\endpspicture
\pspicture(-1,0)(3,2.5)
\pspolygon[fillstyle=solid,fillcolor=lightgray](0,0)(2,0)(1,2)(0,0)
\uput[dl](0,0){\small $a$}
\rput(2.1,-0.1){ \small $b$}
\uput[u](1,2){\small $c$}
\uput[d](1,0){\small $ab$}
\uput[r](1.4,1){\small $bc$}
\uput[r](-0.2,1){\small $ac$}
\rput(2.3,2.6){\small $abc$}
\psline{->}(2.1,2.2)(1,1.1)
\rput(1,-1.1){\large $\Delta_{\mathscr{G}}(K_3)$}
\endpspicture
\pspicture(-3,0)(2,2.5)
\uput[dl](0,0){$a$}
\uput[dr](2,0){ $b$}
\uput[ur](2,2){$c$}
\uput[ul](0,2){$d$}
\rput(0,0){ $\bullet$}
\rput(2,0){ $\bullet$}
\rput(2,2){ $\bullet$}
\rput(0,2){ $\bullet$}
\pspolygon(0,0)(2,0)(2,2)(0,2)
\rput(1,-1.1){\large $C_4$}
\endpspicture
\pspicture(-2,0)(2,2)
\pspolygon(0,0)(2,0)(2,2)(0,2)
\rput(-0.2,-0.1){ \small $a$}
\rput(2.1,-0.1){ \small $b$}
\rput(2.1,2.1){ \small $c$}
\rput(-0.2,2.1){ \small $d$}
\rput(1,-0.2){ \small $ab$}
\rput(1,2.2){ \small $cd$}
\rput(2.2,1){ \small $bc$}
\rput(-0.4,1){ \small $ad$}
\rput(1,-1.1){\large $\Delta_{\mathscr{G}}(C_4)$}
\endpspicture
\vspace{2 cm}

\psset{unit=0.3cm}
\pspicture(-2,-3)(16,5)
\pspolygon(0,0)(12,0)(12,4)(0,4)(0,0)
\psline(0,4)(4,0)(8,4)(12,0)
\psline(0,0)(4,4)(4,0)
\psline(8,4)(8,0)(12,4)
\pscurve(12,4)(6,6)(0,4)
\rput(0,0){$\bullet$}
\rput(0,4){$\bullet$}
\rput(4,0){$\bullet$}
\rput(4,4){$\bullet$}
\rput(8,0){$\bullet$}
\rput(8,4){$\bullet$}
\rput(12,0){$\bullet$}
\rput(12,4){$\bullet$}
\uput[ul](0,4){\small $d$}
\uput[u](4,4){\small $c$}
\uput[u](8,4){\small $e$}
\uput[ur](12,4){\small $h$}
\uput[dl](0,0){\small $a$}
\uput[d](4,0){\small $b$}
\uput[d](8,0){\small $f$}
\uput[dr](12,0){\small $g$}
\rput(6,-3){\large $H$}
\endpspicture
\psset{unit=0.25cm}
\pspicture(-5,-6)(14,4)
\pspolygon[fillstyle=solid,fillcolor=yellow](0,-3)(3,0)(8,0)
\pspolygon[fillstyle=solid,fillcolor=red](0,-3)(11,-3)(8,0)
\pspolygon[fillstyle=solid,fillcolor=lightgray](-3,0)(0,-3)(3,0)(0,4)
\psline(0,-3)(0,4)
\psline[linestyle=dashed](-3,0)(3,0)
\pspolygon[fillstyle=solid,fillcolor=lightgray](8,0)(11,-3)(14,0)(11,4)
\psline(11,-3)(11,4)
\psline[linestyle=dashed](8,0)(14,0)
\psline(0,4)(11,4)
\uput[l](-3,0){\small $a$}
\uput[u](0,4){\small $d$}
\uput[dr](0,-3){\small $b$}
\uput[ur](3,0){\small $c$}
\uput[ul](8,0){\small $e$}
\uput[u](11,4){\small $h$}
\uput[d](11,-3){\small $f$}
\uput[dr](14,0){\small $g$}
\rput(6,-6.3){\large $\Delta_{\mathscr{G}}(H)$}
\endpspicture
  \end{center}
\caption{Illustration of $\Delta_{\ccG}$}
\label{Illustration_Delta_G}
\end{figure}

\begin{lemma}[{\cite[Proposition 3.2]{prisner92}}] \label{lemma-prisner}
Let $G$ be a graph  and  $g$ be a dismantlable vertex in  $G$; 
then  $\Delta_{\mathscr{G}}(G) \red \Delta_{\mathscr{G}}(G\setminus g)$. 
\end{lemma}

\proof
 Let $a$ be a vertex which dominates $g$ (\textit{i.e.} $N_G[g]\subseteq N_G[a]$ 
with $a\neq g$);
for every maximal complete subgraph $c$ of $G$, we have $g \in c \Rightarrow a \in c$.
So, let  $c=[g,a,g_1,\ldots,g_n]$ a maximal complete subgraph of $G$ which contains $g$; 
then $c'=[g,g_1,\ldots,g_n]$ is a free face
of c (taken as a simplex of $\Delta_{\mathscr{G}}(G))$ and 
$\Delta_{\mathscr{G}}(G) \red \Delta_{\mathscr{G}}(G)\setminus \{c,c'\}$.
By iteration, $\Delta_{\mathscr{G}}(G)$ collapses on the subcomplex 
formed by all simplices which do not contain 
$g$, \textit{i.e.} $\Delta_{\mathscr{G}}(G)\red \Delta_{\mathscr{G}}(G\setminus g)$.
\endproof

\begin{proposition}\label{s_prop_G_to_K} Let $G,H \in \mathscr{G}$. 
Then,  $G \red H \Longrightarrow \Delta_{\mathscr{G}}(G) \red \Delta_{\mathscr{G}}(H)$.
\end{proposition}

\proof
It suffices to prove that if $g$ is s-dismantlable in  $G$, then 
$\Delta_{\mathscr{G}}(G) \red \Delta_{\mathscr{G}}(G\setminus g)$.
It follows from Lemma  \ref{lemma-prisner} that $\Delta_{\mathscr{G}}(H)$ is 
collapsible for every dismantlable graph  $H$.
 Thus, by definition of s-dismantlability,  
${\rm link}_{\Delta_{\mathscr{G}}(G)}(<g>)=\Delta_{\mathscr{G}}(N_G(g))$ is collapsible
 when $g$ is s-dismantlable (where $<g>$ is the 0-simplex 
of $\Delta_{\mathscr{G}}(G)$
determined by $g$). As
$\Delta_{\mathscr{G}}(G\setminus g)=\Delta_{\mathscr{G}}(G) 
\setminus {\rm star}^o_{\Delta_{\mathscr{G}}(G)}(<g>)$, 
the  conclusion follows from Lemma \ref{lemma_welker}.
\endproof

\begin{remark} \label{s_converse_G_to_K_false}
 The converse of Proposition \ref{s_prop_G_to_K}
 is not true;  a counterexample is given by the graphs $G$ and $H$ 
of Fig. \ref{ctex_s_Delta_G}.
Indeed, we have $K\red L=K\setminus \{ <a,b,c>,<b,c>\}$ (the complex $K$ 
collapses onto $L$) with 
$K=\Delta_{\mathscr{G}}(G)$, $L=\Delta_{\mathscr{G}}(H)\equiv H$ and we 
don't have $G \red H$ because
there is no s-dismantlable vertex in $G$.

\begin{figure}[h]
\psset{unit=0.5 cm}
\pspicture(-5,-2.5)(2,2.8)
\psline(0,-2)(-2,-2)(-2,0)(2,0)(2,2)(0,2)(0,-2)
\rput(0,0){$\bullet$}
\uput[dr](0,0){$a$}
\rput(2,0){$\bullet$}
\rput(2,2){$\bullet$}
\rput(0,2){$\bullet$}
\rput(0,-2){$\bullet$}
\rput(-2,0){$\bullet$}
\rput(-2,-2){$\bullet$}
\uput[ul](0,2){$b$}
\uput[ul](-2,0){$c$}
\pspolygon[fillstyle=solid,fillcolor=lightgray](0,0)(0,2)(-2,0)
\rput(1.5,-1.5){$K$}
\endpspicture
\pspicture(-4,-2.5)(2,2.5)
\psline(0,-2)(-2,-2)(-2,0)(2,0)(2,2)(0,2)(0,-2)
\rput(0,0){$\bullet$}
\uput[dr](0,0){$a$}
\rput(2,0){$\bullet$}
\rput(2,2){$\bullet$}
\rput(0,2){$\bullet$}
\rput(0,-2){$\bullet$}
\rput(-2,0){$\bullet$}
\rput(-2,-2){$\bullet$}
\uput[ul](0,2){$b$}
\uput[ul](-2,0){$c$}
\rput(1.5,-1.5){$L$}
\endpspicture
\pspicture(-5,-2.5)(2,2.5)
\psline(0,-2)(-2,-2)(-2,0)(2,0)(2,2)(0,2)(0,-2)
\rput(0,0){$\bullet$}
\uput[dr](0,0){$a$}
\rput(2,0){$\bullet$}
\rput(2,2){$\bullet$}
\rput(0,2){$\bullet$}
\rput(0,-2){$\bullet$}
\rput(-2,0){$\bullet$}
\rput(-2,-2){$\bullet$}
\psline(0,2)(-2,0)
\rput(1.5,-1.5){$G$}
\uput[ul](0,2){$b$}
\uput[ul](-2,0){$c$}
\endpspicture
\pspicture(-4,-2.5)(2,2.5)
\psline(0,-2)(-2,-2)(-2,0)(2,0)(2,2)(0,2)(0,-2)
\rput(0,0){$\bullet$}
\uput[dr](0,0){$a$}
\rput(2,0){$\bullet$}
\rput(2,2){$\bullet$}
\rput(0,2){$\bullet$}
\rput(0,-2){$\bullet$}
\rput(-2,0){$\bullet$}
\rput(-2,-2){$\bullet$}
\uput[ul](0,2){$b$}
\uput[ul](-2,0){$c$}
\rput(1.5,-1.5){$H$}
\endpspicture
\caption{$\Delta_{\mathscr{G}}(G) \red \Delta_{\mathscr{G}}(H)$ but $G \not\red H$}
\label{ctex_s_Delta_G}
\end{figure}
\end{remark}

\subsection{From $\ccK$ to $\ccG$}

We consider the application  $\Gamma : \ccK \to \mathscr{G}$ whose definition is :
if $K \in \ccK$, $\Gamma(K)$ is the graph whose vertices are the simplices of $K$
with  edges $\{\s,\s'\}$ when $\s \subset \s'$ or  $\s' \subset \s$.
\ss

\noindent If $\sigma$ is a simplex of $K \in \ccK$, we shall write $K[\sigma]$ 
for the simplicial
subcomplex of $K$ formed by all faces of 
$\sigma$ ($K[\sigma]:=\{\tau \in K\:,\: \tau \subset \sigma\}$);
if $\tau$ is a maximal face of $\sigma$, $K[\sigma]\setminus\{\sigma,\tau\}$ 
is a simplicial complex.
In order to understand the  relation of formal deformations to s-dismantlability, 
we have the following results:

\begin{lemma}\label{lemma_K_to_G} If $\tau$ is a maximal face of  $\sigma$ then
$\Gamma(K[\sigma]\setminus\{\sigma,\tau\})$ is dismantlable.
\end{lemma}

\proof
 Let $\sigma=<a_0,a_1,\ldots,a_n>$ and $\tau=<a_1,\ldots,a_n>$. 
The vertices of $\Gamma(K[\sigma]\setminus\{\sigma,\tau\})$ are the simplices of
$K[\sigma]\setminus\{\sigma,\tau\}$.
These vertices can be written $<a_{i_1},\ldots,a_{i_k}>$
with $0\leq i_1 < i_2<\ldots <i_k\leq n$
excepting $<a_0,a_1,\ldots,a_n>$ and $<a_1,\ldots,a_n>$;
we shall say that such a vertex $<a_{i_1},\ldots,a_{i_k}>$  contains  $a$ 
if $a \in \{a_{i_1},\ldots,a_{i_k}\}$.
Every vertex  $ x=<a_{i_1},\ldots,a_{i_{n-1}}>$ which does 
not contain $a_0$ is dismantlable
(because it is dominated by $<a_0,a_{i_1},\ldots,a_{i_{n-1}}>$
the unique $(n-1)$-simplex containing  $<a_{i_1},\ldots,a_{i_{n-1}}>$).
Thus, $\Gamma(K[\sigma] \setminus \{\sigma,\tau\})$ can be 
dismantled on $\Gamma_{n-1}$ obtained by deleting
all vertices corresponding to  $(n-1)$-simplices which do not contain  $a_0$.
Next, every vertex  $<a_{i_1},\ldots,a_{i_{n-2}}>$ which does 
not contain $a_0$ is dismantlable in $\Gamma_{n-1}$
(because it is dominated by $<a_0,a_{i_1},\ldots,a_{i_{n-2}}>$
the unique $(n-2)$-simplex containing  $<a_{i_1},\ldots,a_{i_{n-2}}>$).
Thus, $\Gamma_{n-1}$ can be dismantled on $\Gamma_{n-2}$ obtained by deleting
all vertices corresponding to  $(n-2)$-simplices which do not contain  $a_0$.
The iteration of this procedure shows that 
$\Gamma(K[\sigma] \setminus \{\sigma,\tau\})$ 
is dismantlable on its subgraph induced by the vertices containing $a_0$. 
But this subgraph is a cone
on $<a_0>$ and  this shows that  $\Gamma(K[\sigma] \setminus \{\sigma,\tau\})$  
is dismantlable.
\endproof

\begin{proposition}\label{prop_K_to_G} 
Let $K,L \in \ccK$. Then, $K \red L \Longrightarrow \Gamma(K) \red \Gamma(L)$.
\end{proposition}

\proof
It suffices to prove that if $\{\sigma,\tau\}$ is a collapsible pair in $K$, then 
$\Gamma(K)\red \Gamma(K\setminus \{\sigma,\tau\})$. We note that 
$ \Gamma(K\setminus \{\sigma,\tau\})= \Gamma(K) \setminus \{\sigma,\tau\}$
and  that the vertex $\tau$ is dismantlable in  $\Gamma(K)$ 
(because it is dominated by $\sigma$);
so, we have the reduction
$\Gamma(K) \red \Gamma(K) \setminus \tau$.
Now, we have  $N_{\Gamma(K)\setminus \tau}(\sigma)=
\Gamma(K[\sigma]\setminus \{\sigma,\tau\})$
(because $\sigma$ is a maximal simplex),
and we conclude that $\sigma$ is s-dismantlable by the Lemma \ref{lemma_K_to_G}.
\endproof

\subsection{Barycentric subdivision}

Let us recall the notion of barycentric subdivision in $\ccK$.  If  $K \in \ccK$,  
the $n$-simplices of the barycentric subdivision  $Bd(K)$ (or $K'$) of $K$
are the  $<\s_0,\s_1,\ldots,\s_n>$ composed of  $n+1$
simplices of $K$ such that $\s_0 \subset \s_1 \subset \ldots
\subset \s_n$.
Now, we define a similar notion in $\ccG$.

\begin{definition}
If $G \in \mathscr{G}$, the barycentric subdivision  $Bd(G)$ (or $G'$) of $G$ is the graph
whose vertices are the complete subgraphs of $G$ and there is an edge 
between two vertices if, and only if, there is
an inclusion between the two corresponding complete subgraphs.
\end{definition}

\begin{remark} Each complete subgraph  of cardinality at least two creates a new vertex 
in the barycentric subdivision (cf. Fig. \ref{dessin_barycentric_subdiv}).
The equalities $\Gamma \circ \Delta_{\mathscr{G}} = Bd$ (in $\mathscr{G}$)
and $\Delta_{\mathscr{G}} \circ \Gamma = Bd$ (in $\ccK$) follow directly 
from the definitions and will be useful in the  following.
\end{remark}

\begin{figure}[h]
\psset{unit=0.6 cm}
\pspicture(-9,-2)(4,3.5)
\rput(-5,0){\large $\bullet$}
\rput(-2,0){\large $\bullet$}
\rput(0,2){\large $\bullet$}
\rput(2,0){\large $\bullet$}
\rput(-5,2){\large $\bullet$}
\rput(-7,2){\large $\bullet$}
\rput(-7,0){\large $\bullet$}
\psline(-5,2)(-5,0)(-7,0)(-7,2)(-5,2)(-2,0)(0,2)(2,0)(-2,0)
\uput[dl](-7,0){$a$}
\uput[ul](-7,2){$b$}
\uput[ur](-5,2){$c$}
\uput[dr](-5,0){$d$}
\uput[dl](-2,0){$e$}
\uput[u](0,2){$f$}
\uput[dr](2,0){$g$}
\rput(-3.5,-1.5){\large $G$}
\endpspicture
\pspicture(-9,-2)(4,3.5)
\rput(-5,0){\large $\bullet$}
\rput(-2,0){\large $\bullet$}
\rput(0,2){\large $\bullet$}
\rput(2,0){\large $\bullet$}
\rput(-5,2){\large $\bullet$}
\rput(-7,2){\large $\bullet$}
\rput(-7,0){\large $\bullet$}
\psline(-5,2)(-5,0)(-7,0)(-7,2)(-5,2)(-2,0)(0,2)(2,0)(-2,0)
\uput[dl](-7,0){$a$}
\uput[ul](-7,2){$b$}
\uput[ur](-5,2){$c$}
\uput[dr](-5,0){$d$}
\uput[l](-7,1){$ab$}
\uput[r](-5,1){$cd$}
\uput[u](-6,2){$bc$}
\uput[d](-6,0){$ad$}
\uput[ur](-3.5,1){$ce$}
\uput[dl](-2,0){$e$}
\uput[u](0,2){$f$}
\uput[dr](2,0){$g$}
\rput(-6,0){\large $\bullet$}
\rput(-7,1){\large $\bullet$}
\rput(-6,2){\large $\bullet$}
\rput(-5,1){\large $\bullet$}
\rput(-3.5,1){\large $\bullet$}
\rput(-1,1){\large $\bullet$}
\rput(1,1){\large $\bullet$}
\rput(0,0){\large $\bullet$}
\rput(0,0.8){\large $\bullet$}
\uput[ur](1,1){$fg$}
\uput[ul](-1,1){$ef$}
\uput[d](0,0){$eg$}
\uput[ur](2.2,2.7){$efg$}
\psline[linestyle=dashed,linecolor=lightgray]{->}(2.2,2.7)(0.1,0.9)
\psline(-2,0)(0,0.8)(-1,1)
\psline(2,0)(0,0.8)(1,1)
\psline(0,0)(0,2)
\rput(-3.5,-1.5){\large $G'=Sd(G)$}
\endpspicture
\caption{Example of barycentric subdivision}
\label{dessin_barycentric_subdiv}
\end{figure}

\begin{proposition} \label{prop_subdivision} For every $G \in \mathscr{G}$,  
$G$ and $G'$ have the same s-homotopy type (\textit{i.e.} $[G]_s=[G']_s$).
\end{proposition}

\proof
Let $n$ be the  cardinal of $V(G)$; we choose to number the vertices of $G$;
 thus, we have $V(G)=\{g_1,g_2,\ldots,g_n\}$. 
Let us recall that   $V(G')=\ccC(G)$, the set of  complete subgraphs of $G$.
 In what follows, every complete subgraph $c$ is considered under its unique expression
$c=[g_{i_1},\ldots,g_{i_k}]$ with $1\leq i_1<i_2<\ldots <i_k\leq n$; 
we shall denote by $i_k ={\rm max\:}(c)$.
If $g \in \{g_{i_1},\ldots,g_{i_k}\}$ and $c =[g_{i_1},\ldots,g_{i_k}]\in \ccC(G)$,
we shall write $g \in c$ and for  $c,d \in \ccC(G)$, we shall write $c \subset d$ when
$c =[g_{i_1},\ldots,g_{i_k}]$, $d =[g_{j_1},\ldots,g_{j_m}]$ 
and $\{i_1,\ldots,i_k\} \subset \{j_1,\ldots,j_m\}$. 
We know that 
$E(G')=\{\: cd\:,\: c,d \in \ccC(G)~{\rm with}~c \subset d~{\rm or}~d \subset c\}$.
To prove $[G]_s=[G']_s$, we go from  $G$ to $G'$ in two steps 
(addition and  suppression of s-dismantlable vertices).
\ss

First step: For every $c \in \ccC(G)$, we add a vertex $\widehat c$ to $G$. 
We begin with complete subgraphs of cardinal 1, 
we proceed with complete subgraphs of cardinal 2, 
next with complete subgraphs of cardinal 3... 
until we have reached all complete subgraphs. When we add a vertex $\widehat c$ 
corresponding to the complete subgraph 
 $c=[g_{i_1},\ldots,g_{i_k}]$  of cardinal $k$, we add the edges 
 $\widehat c \widehat d$ if $d\subset c$, $\widehat c g_{i_k}$ and $\widehat c  g_j$ if
$j >i_k$ and  $c \cup g_j \in \ccC(G)$; 
this corresponds to the addition of an s-dismantlable vertex
because the open neighborhood of
 $\widehat c$ (when we add it) is a cone on $g_{i_k}$.
\ss

The  graph $H$ obtained at the end of the first step  is such that $V(H)=V(G)\cup V(G')$.
\ss

Second step:  We note that $g_1$ is s-dismantlable in $H$
(because $N_H(g_1)=\{\widehat{[g_1]}\} \cup N_G(g_1)$ is a cone on $\widehat{[g_1]}$)
and, more generally, for $2 \leq i \leq n$, 
let us verify that the vertex  $g_i$ is s-dismantlable 
in $H_i =H \setminus \{g_1,\ldots,g_{i-1}\}$. We have $W_i=N_{H_i}(g_i)=
\{ \widehat c\:, \: c\in \ccC(G), c\cup g_i \in \ccC(G) ~
{\rm and~} \: {\rm max\:} (c) \leq  i\} \cup \{g_j\:,\:
g_j \in N_G(g_i)~{\rm and}~j>i\}$. Let us denote $W'_i=
\{ \widehat c\:, c \in \ccC(G) ~{\rm and~} \: {\rm max\:} (c) =  i\} \cup \{g_j\:,\:
g_j \in N_G(g_i)~{\rm and}~j>i\}$, the cone on $\widehat{[g_i]}$. We have either $W_i=W'_i$, or $W_i$ 
is dismantlable on $W'_i$. Indeed, let us suppose $W_i\neq W'_i$ 
 and let $\widehat c$ 
such that $c\cup g_i \in \ccC(G) ~{\rm and~} \: {\rm max\:} (c) <  i$.
Then $\widehat{c\cup g_i} \in W_i$ and $N_{W_i}(\widehat{c}) 
\subseteq N_{W_i}(\widehat{c\cup g_i})$; so, $\widehat c$ is dismantlable in $W_i$
and more generally $W_i$ is dismantlable on $W'_i$.
Consequently, $W_i$ is dismantlable, \textit{i.e.} $g_i$ is s-dismantlable in $H_i$.
Thus, in $H$, one can s-delete all vertices of $G$ (in the following order: $g_1$, $g_2$,
$\ldots$, $g_n$); the resultant graph is $G'$.
\endproof

\subsection{Correspondence of homotopy classes}

\begin{theorem}\label{theoremGK} 1. 
Let $G,H \in \mathscr{G}$; $G$ and $H$ have the same s-homotopy type
if, and only if, $\Delta_{\mathscr{G}}(G)$ and $\Delta_{\mathscr{G}}(H)$ 
have the same simple-homotopy type:
$$[G]_s =[H]_s \Longleftrightarrow [\Delta_{\mathscr{G}}(G)]_s 
=[\Delta_{\mathscr{G}}(H)]_s$$

2.  Let $K,L \in \ccK$; $K$ and $L$ have the same simple-homotopy type
if, and only if, $\Gamma(K)$ and $\Gamma(L)$ have the same s-homotopy type:
$$[K]_s =[L]_s \Longleftrightarrow [\Gamma(K)]_s =[\Gamma(L)]_s$$
\end{theorem}

\proof
 1. \udl{$\Longrightarrow$ :} corollary of Proposition \ref{s_prop_G_to_K}.
\udl{$\Longleftarrow$ :} By the Proposition \ref{prop_K_to_G}, we get 
$[\Delta_{\mathscr{G}}(G)]_s =[\Delta_{\mathscr{G}}(H)]_s$ $\Longrightarrow$ 
$[G']_s=[\Gamma(\Delta_{\mathscr{G}}(G))]_s 
=[\Gamma(\Delta_{\mathscr{G}}(H))]_s=[H']_s$
and we conclude 
with the Proposition \ref{prop_subdivision}.

2. \udl{$\Longrightarrow$ :} corollary of Proposition \ref{prop_K_to_G}.
\udl{$\Longleftarrow$ :} By using assertion 1 of the theorem, we obtain
$ [\Gamma(K)]_s=[\Gamma(L)]_s \Longrightarrow 
[\Delta_{\mathscr{G}}(\Gamma(K))]_s=[\Delta_{\mathscr{G}}(\Gamma(L))]_s$. 
So, we get $[K']_s=[L']_s$ and we can conclude  $[K]_s = [L]_s$ 
because it is well known that a simplicial complex and its barycentric subdivision 
have the same simple-homotopy type (\cite{cohen}).
\endproof

\begin{remark}
It is  clear that an s-collapsible graph is s-dismantlable 
(\textit{i.e.}, $G \red pt \Rightarrow [G]_s=[pt]_s$) and 
it follows from Proposition \ref{s_prop_G_to_K} and Theorem \ref{theoremGK} 
that there exists s-dismantlable and not s-collapsible graphs. 
Indeed, there are well known examples of simplicial complexes $K$
(triangulations of the dunce hat (\cite{zeeman}) or of the 
Bing's house (\cite{cyy}), for instance) 
which are not collapsible ($K \not\red pt$) but have the same simple-homotopy 
type of a point ($[K]_s=[pt]_s$).
So, for example, the  graph $D$ of Fig. \ref{dunce_hat} 
(with 17 vertices and 36 triangles)  is s-dismantlable but not s-collapsible 
because $\Delta_{\ccG}(D)$ is a triangulation of the dunce hat. 
This graph $D$ is actually the comparability
 graph of a poset given in \cite[Figure 2, p.380]{walk81} and named here $P_d$.
An example of a graph $B$ such that $\Delta_{\ccG}(B)$ is a triangulation 
of the Bing's house is given 
in \cite[\S 5]{cyy}.

\begin{figure}[h]
\psset{unit=0.2 cm}
\pspicture(-20,-1.5)(13,20)
\pspolygon(-12,0)(12,0)(0,18)(-12,0)
\pspolygon(2,11)(6,5)(4,2)(-4,2)(-6,5)(-2,11)(2,11)
\rput(-12,0){$\bullet$}
\rput(12,0){$\bullet$}
\rput(-4,0){$\bullet$}
\rput(0,0){$\bullet$}
\rput(4,0){$\bullet$}
\rput(8,6){$\bullet$}
\rput(6,9){$\bullet$}
\rput(4,12){$\bullet$}
\rput(0,18){$\bullet$}
\rput(-8,6){$\bullet$}
\rput(-6,9){$\bullet$}
\rput(-4,12){$\bullet$}
\rput(2,11){$\bullet$}
\rput(4,8){$\bullet$}
\rput(6,5){$\bullet$}
\rput(-2,11){$\bullet$}
\rput(-4,8){$\bullet$}
\rput(-6,5){$\bullet$}
\rput(-4,2){$\bullet$}
\rput(0,2){$\bullet$}
\rput(4,2){$\bullet$}
\rput(5,3.5){$\bullet$}
\rput(-5,3.5){$\bullet$}
\rput(0,11){$\bullet$}
\rput(0,6){$\bullet$}
\psline(0,0)(0,18)\psline(-12,0)(6,9)\psline(12,0)(-6,9)
\psline(0,6)(2,11)\psline(0,6)(6,5)
\psline(0,6)(-2,11)\psline(0,6)(-6,5)
\psline(0,6)(-4,2)\psline(0,6)(4,2)
\psline(8,6)(6,5)(6,9)(2,11)(4,12)\psline(-8,6)(-6,5)(-6,9)(-2,11)(-4,12)\psline(-4,0)(-4,2)(0,0)(4,2)(4,0)
\psline(4,2)(12,0)(6,5)\psline(-4,2)(-12,0)(-6,5)\psline(2,11)(0,18)(-2,11)
\uput[dl](-12,0){1}\uput[dr](12,0){1}\uput[u](0,18){1}
\uput[d](-4,0){2}\uput[ur](4,12){2}\uput[ul](-8,6){2}
\uput[d](4,0){3}\uput[ul](-4,12){3}\uput[ur](8,6){3}
\uput[d](0,0){4}\uput[ur](6,9){4}\uput[ul](-6,9){4}
\rput(-12,16){\large $D$}
\endpspicture
\pspicture(-28,-1.5)(13,20)

\rput(-12,0){$\bullet$}
\uput[d](-12,0){4}
\rput(12,0){$\bullet$}
\uput[d](12,0){1}
\rput(0,0){$\bullet$}

\rput(-12,9){$\bullet$}
\rput(-9,9){$\bullet$}
\rput(-15,9){$\bullet$}
\rput(12,9){$\bullet$}
\rput(9,9){$\bullet$}
\rput(15,9){$\bullet$}
\rput(-2,9){$\bullet$}
\uput[l](-2,9){2}
\rput(2,9){$\bullet$}
\uput[r](2,9){3}

\rput(-12,18){$\bullet$}
\rput(-9,18){$\bullet$}
\rput(-6,18){$\bullet$}
\rput(12,18){$\bullet$}
\rput(9,18){$\bullet$}
\rput(6,18){$\bullet$}

\psline(-12,0)(-15,9)\psline(-12,0)(-12,9)\psline(-12,0)(-9,9)\psline(-12,0)(-2,9)\psline(-12,0)(2,9)
\psline(12,0)(15,9)\psline(12,0)(12,9)\psline(12,0)(9,9)\psline(12,0)(-2,9)\psline(12,0)(2,9)
\psline(0,0)(-12,9)\psline(0,0)(-15,9)\psline(0,0)(-9,9)
\psline(0,0)(12,9)\psline(0,0)(15,9)\psline(0,0)(9,9)
\psline(-2,9)(-12,18)\psline(-2,9)(-9,18)\psline(-2,9)(12,18)
\psline(2,9)(-6,18)\psline(2,9)(6,18)\psline(2,9)(9,18)
\psline(-12,18)(-15,9)(6,18)(15,9)(-6,18)(-9,9)(12,18)(12,9)
(-9,18)(-12,9)(9,18)(9,9)(-12,18)

\rput(-17,16){\large $P_d$}
\endpspicture
\caption{$\Delta_{\ccG}(D)$ is a triangulation of the dunce hat: $D=Comp(P_d)$}
\label{dunce_hat}
\end{figure}

\end{remark}

\section{Relation with posets}\label{sec_poset}

\subsection{From $\ccP$ to $\ccG$}

Let  $\ccP$ be the set of finite partially ordered sets or finite posets. 
In what follows, when $P\subseteq Q$ with $Q \in \ccP$,
$P$ will be called {\sl subposet of $Q$} if, 
for every $x,y$ in $P$, $x\leq_Py \Longleftrightarrow x\leq_Q y$.
If $P \in \ccP$, $Comp(P)\in \mathscr{G} $ is the comparability graph of $P$ 
(its vertices are the elements of $P$ with an edge $xy$ if, and only if, 
$x$ and $y$ are comparable).

Let $P\in \ccP$. For every $x$ in $P$, we define
$P_{<x}:=\{y\in P\:,\: y<x\}$ and $P_{>x}:=\{y\in P\:,\: y>x\}$.
We recall that $x$ is  {\sl irreducible}\footnote{It seems to be the most 
classical terminology 
(\cite{rival},\cite{bcf94},\cite{ginsburg}); in 
 \cite{barmin}, irreducible points are called {\sl (up {\rm or} down) beat points}.} 
either if $P_{<x}$ 
has a maximum, or if
$P_{>x}$ has a minimum. The poset $P$ is called {\sl dismantlable} 
if we can write $P=\{x_1,x_2,\ldots,x_n\}$
 with $x_i$ irreducible in the subposet induced 
by $\{x_1,\ldots,x_i\}$, for $2 \leq i \leq n$.
Let us recall that a {\sl cone} is a poset having a maximum or a minimum;
 if we can write $P=P_{\geq x}$ or
$P=P_{\leq x}$ for some $x$ in $P$, $P$ will be called a cone \emph{on} $x$. 
Cones are examples of dismantlable posets.

If $P,Q \in \ccP$, $P*Q$ is the poset whose elements are those of $P$ and $Q$ and with
the relations $p\leq _P p'$, $q\leq _Q q'$ and $p\leq q$ for all $p,p'\in P$ 
and $q,q'\in Q$.
In particular, $P*\emptyset=\emptyset * P=P$ for all $P\in \ccP$.

\begin{lemma} \label{lemma_star}
Let $P,Q\in \ccP$;  $P * Q$ is dismantlable if, and only if, $P$ or $Q$ is dismantlable.
\end{lemma}

\proof
Let us suppose that $P*Q$ is dismantlable with
$P*Q=\{x_1,x_2,\ldots,x_N\}$ (where $N=\vert P\vert + \vert Q \vert$
and $x_i$ is irreducible in the subposet of $P*Q$ induced by $\{x_1,\ldots,x_i\}$, 
for $2\leq i \leq N$)
and that $Q$ is a non dismantlable poset. 
We can write $P=\{x_{i_1},x_{i_2},\ldots,x_{i_k}\}$
(where $k=\vert P\vert$) with $i_j<i_l$ for all $1\leq j<l\leq k$.
We shall verify that $P$ is dismantlable with $x_{i_l}$ irreducible in the subposet of $P$ 
induced by $P\setminus \{x_{i_{l+1}},\ldots,x_{i_k}\}=\{x_{i_1},\ldots,x_{i_l}\}$, 
for $2\leq l \leq k$.
We know that $x_{i_l}$ is irreducible 
in $(P*Q)\setminus \{x_{i_l+1},x_{i_l+2},\ldots,x_{N}\}$.
\begin{itemize}
\item First case: a maximum of
$\bigl((P*Q)\setminus \{x_{i_l+1},x_{i_l+2},\ldots,x_{N}\}\bigr)_{<x_{i_l}}$ is 
also a maximum of
$\bigl(P\setminus \{x_{i_{l+1}},x_{i_{l+2}},\ldots,x_{i_k}\}\bigr)_{<x_{i_l}}$ 
(this follows from
$\bigl((P*Q)\setminus \{x_{i_l+1},x_{i_l+2},\ldots,x_{N}\}\bigr)_{<x_{i_l}}=
\bigl(P\setminus \{x_{i_{l+1}},x_{i_{l+2}},\ldots,x_{i_k}\}\bigr)_{<x_{i_l}}$). 
\item Second case: a minimum of
$\bigl((P*Q)\setminus \{x_{i_l+1},x_{i_l+2},\ldots,x_{N}\}\bigr)_{>x_{i_l}}$ 
is also a minimum of
$\bigl(P\setminus \{x_{i_{l+1}},x_{i_{l+1}},\ldots,x_{i_k}\}\bigr)_{>x_{i_l}}$. 
This follows from 
$\bigl((P*Q)\setminus \{x_{i_l+1},x_{i_l+2},\ldots,x_{N}\}\bigr)_{>x_{i_l}}=
\bigl(P\setminus \{x_{i_{l+1}},x_{i_{l+2}},\ldots,x_{i_k}\}\bigr)_{>x_{i_l}}*Q'$
where $Q'$ is a subposet of $Q$ and the fact  that
$\bigl(P\setminus \{x_{i_{l+1}},x_{i_{l+1}},\ldots,x_{i_k}\}\bigr)_{>x_{i_l}} 
\neq \emptyset$.
Indeed, $\bigl((P*Q)\setminus \{x_{i_l+1},x_{i_l+2},\ldots,x_{N}\}\bigr)_{>x_{i_l}}=Q'$
would mean that $Q'$ is dismantlable (because
$\bigl((P*Q)\setminus \{x_{i_l+1},x_{i_l+2},\ldots,x_{N}\}\bigr)_{>x_{i_l}}$ 
has a minimum and, thus, is dismantlable)
which contradicts the fact that $Q$ is non dismantlable (because $Q'$
 is a subposet of $Q$ obtained by suppression of irreducible elements).
\end{itemize}
In conclusion, $x_{i_l}$ is irreducible in $\{x_{i_1},\ldots,x_{i_l}\}$, 
for $2\leq l \leq k$ and $P$ is dismantlable.

Reciprocally, let us suppose that $P$ is a dismantlable poset and that 
we have $P=\{x_1,x_2,\ldots,x_n\}$
 with $x_i$ irreducible in the subposet induced 
by $\{x_1,\ldots,x_i\}$, for $2\leq i \leq n$. 
By the equalities $(P*Q)_{>p}=P_{>p}*Q$ and $(P*Q)_{<p}=P_{<p}$,
we see that $x_i$ is irreducible in the subposet of $P*Q$ induced 
by $\{x_1,\ldots,x_i\}*Q$ if $x_i$ is irreducible in the subposet of $P$ induced 
by $\{x_1,\ldots,x_i\}$. As a consequence, we can go from $P*Q$ to $\{x_1\}*Q$ 
by successive suppressions
of irreducibles and this shows that $P*Q$ is dismantlable (because $\{x_1\}*Q$ 
is a cone). A similar argument
yields that $P*Q$ is dismantlable if we suppose $Q$ dismantlable.
\endproof

In \cite{barmin}, Barmak and Minian introduce the notion of {\sl weak points} in a poset:
 $x\in P$ is a {\sl weak point} if
$P_{<x}$ or $P_{>x}$ is dismantlable.
So, by Lemma \ref{lemma_star}, $x$ is a weak point if, and only if, $P_{>x} *P_{<x}$ 
is dismantlable.
Now, it is well known  (\cite{bcf94},\cite{ginsburg}) that a poset $P$ is a 
dismantlable poset if, and only if, 
$Comp(P)$ is a dismantlable graph. 
As we have $N_{Comp(P)}(x)=Comp(P_{>x} * P_{<x})$, we obtain the 
following result:

\begin{proposition}  \label{prop_weak_and_s-dismant}
Let $P\in \ccP$ and $x\in P$. Then, 
 $x$ is a {\sl weak point} of $P$  if, and only if, $x$ is s-dismantlable in $Comp(P)$.
\end{proposition}

The notation $P \red P\setminus \{x\}$ will mean that $x$ is a weak point of $P$ 
and we shall write $P \red Q$ if $Q$ is a subposet of $P$ 
obtained by successive deletions of weak points.

\subsection{From $\ccG$ to $\ccP$}

If $G \in \mathscr{G}$, $C(G)\in \ccP$ is the poset whose elements are the 
complete subgraphs of $G$
ordered by inclusion. Before establishing  the relation between reduction 
by s-dismantlable vertices in $\ccG$ and deletion of weak points in $\ccP$, 
we recall that the poset product $P\times Q$ of two posets 
$P$ and $Q$ is the set $P \times Q$ ordered by
$(p,q)\leq (p',q')$ if $p\leq_Pp'$ and $q\leq _Q q'$ 
for all $(p,q),(p',q') \in P \times Q$. In particular,
$P \times \{a,b,a<b\} $ is the poset  formed by two copies of $P$ 
(namely, $P_a:=P\times \{a\}=\{(p,a),p\in P\}$
and $P_b:=P\times \{b\}=\{(p,b),p\in P\}$) with relations of $P$ 
in the two copies $P_a$ and $P_b$ and the additional relations
$(p,a)\leq (p',b)$ if $p\leq_Pp'$.

\begin{lemma} \label{lemme_graphe_vers_poset}
Let $P\in \ccP$
and $S\in \ccP$  such that  $S$ contains 
$W:=P \times \{a,b,a<b\} $ as a subposet with the two following properties:
$$\forall p \in P\:, ~~~~~~~i)~S_{<(p,b)}=W_{<(p,b)}
~~~~~~~~{\rm and} ~~~~~~~~
ii)~S_{>(p,b)}=W_{>(p,b)}$$
Let $Q$ be the poset obtained from $S$ by adding an element $x$ (not in $S$) 
with the only relations $x<(p,b)$ for all $p\in P$.

If $P$ is a dismantlable poset, then $Q \red Q \setminus \{x,(p,b)\:;\:p\in P\} $.
\end{lemma}

\proof
By definition of $Q$, we have $x\leq_Q y$ 
if, and only if, $y \in P_b=\{(p,b),p\in P\}$. Of course, $P_b$ is isomorphic to $P$
and, if we suppose that $P$ is dismantlable, 
this means that $x$ is a weak point in $Q$, \textit{i.e.}
$Q \red S=Q\setminus \{x\}$. 

Now, let $p$ be an irreducible  element in $P$; we shall verify that $(p,b)$ 
is a weak point in $S$.

\begin{itemize}
\item First case:  we suppose that $P_{<p}$ has a maximum element $M$. We get 
$$\begin{array}{rcl}
S_{<(p,b)}\stackrel{i)}{=}W_{<(p,b)}&=&
\{(p',b), p'<p\}\bigcup \{(p',a), p'\leq p\}\\
&=&\{(p',b), p'<p\}\bigcup \{(p',a), p'\leq M\} \bigcup \{(p,a)\}
\end{array}$$
In $S_{<(p,b)}$, we have $y < (p,a) \Leftrightarrow y\leq (M,a)$; in other words, $(M,a)$ 
is a maximum element of $\left(S_{<(p,b)}\right)_{<(p,a)}$ 
and this shows that $(p,a)$ is an irreducible point
in $S_{<(p,b)}$. Now,  
$S_{<(p,b)}\setminus\{(p,a)\}=\{(p',b), p'<p\}\bigcup \{(p',a), p'\leq M\}$
is a cone on $(M,b)$ (because $y\leq (M,b)$ for all $y$ in $S_{<(p,b)}\setminus\{(p,a)\}$)
and we can conclude that $S_{<(p,b)}$ is a dismantlable poset. 

\item Second case: we suppose that $P_{>p}$ has a minimum element $m$, 
then $S_{>(p,b)}\stackrel{ii)}{=}W_{>(p,b)}$
which is a poset isomorphic to $P_{>p}$; so, it is dismantlable (because it is a cone).
\end{itemize}

\noindent The conclusion of the two cases is that $(p,b)$ is a weak point in $S$;
so, we have $S \red S\setminus \{(p,b)\}$.

Now, let us suppose that $P$ is dismantlable with $P=\{p_1,p_2,\ldots,p_n\}$
 with $p_i$ irreducible in the subposet induced 
by $\{p_1,\ldots,p_i\}$, for $2\leq i \leq n$. 
By iterating the preceding discussion we get
$$
\begin{array}{rcl}Q ~ \red  ~ Q\setminus \{x\} ~ \red  ~ Q\setminus \{x,(p_n,b)\} 
& \red  & Q\setminus \{x,(p_n,b),(p_{n-1},b)\}\\
&\red & \ldots \red Q\setminus \{x,(p_n,b),(p_{n-1},b),\ldots,(p_2,b)\}
\end{array}
$$
By condition $i)$, we see that $y < (p_1,b)$ 
in $Q\setminus \{x,(p_n,b),(p_{n-1},b),\ldots,(p_2,b)\}$
if, and only if, $y \leq (p_1,a)$; so, $(p_1,b)$ is a weak point 
in $Q\setminus \{x,(p_n,b),(p_{n-1},b),\ldots,(p_2,b)\}$
(in fact, it is irreducible) and we have proved 
$$ Q ~ \red ~Q\setminus \{x,(p_n,b),\ldots,(p_2,b)\}
~ \red  ~ Q\setminus \{x,(p_n,b),\ldots,(p_2,b),(p_1,b)\}.
$$
\endproof

\begin{proposition}\label{prop_C_and_s-dismant}
If $g$ is s-dismantlable in $G$, then $C(G) \red C(G\setminus g)$.
\end{proposition}

\proof We apply Lemma \ref{lemme_graphe_vers_poset} with $Q=C(G)$ 
and $P=C(N_G(g))$.
More precisely, with the notations of this lemma, 
we have $x=[g]$, $S=C(G)\setminus\{[g]\}$ and
$W= C(N_G[g])\setminus \{[g]\}$. If we denote by $g_i$ the elements of $N_G(g)$,
the isomorphism between $W= C(N_G[g])\setminus \{[g]\}$ 
and $P \times \{a,b,a<b\} =C(N_G(g)) \times \{a,b,a<b\} $
is given by identifying $([g_1,\ldots,g_n],a)\in P \times \{a,b,a<b\}$ 
with $[g_1,\ldots,g_n] \in W$
and $([g_1,\ldots,g_n],b)\in P \times \{a,b,a<b\}$ with $[g,g_1,\ldots,g_n] \in W$.
Conditions $i)$ and $ii)$ are clearly verified.

By supposing  $g$ s-dismantlable in $G$, we get that $N_G(g)$ is a 
dismantlable graph and $P=C(N_G(g))$ a dismantlable
poset by  \cite[Lemma 2.2]{ginsburg}. So, by Lemma \ref{lemme_graphe_vers_poset}, 
we obtain
$$C(G)\red C(G)\setminus \bigl(\{[g]\}\cup 
\{[g,g_1,\ldots,g_n], [g_1,\ldots,g_n]\in C(N_G(g))\}\bigr)=C(G\setminus g)$$
\endproof

\subsection{Correspondence between s-homotopy and simple equivalence}

 A poset $P$ is said (\cite[Definition 3.4]{barmin}) {\sl simply equivalent} to the poset
$Q$ if we can transform $P$ to $Q$ by a finite sequence of additions 
or deletions of weak  points.
 We denote by $[P]_s$ the equivalence class  of $P$ for this 
relation (and call it the \emph{simple type} of $P$).

\begin{theorem}\label{theoremGP}
1. Let $P,Q \in \mathscr{P}$; 
$[P]_s=[Q]_s ~({\rm in}~\ccP) \Longleftrightarrow [Comp(P)]_s=[Comp(Q)]_s ~({\rm in}~\ccG)$.

2. Let  $G,H \in \mathscr{G}$;
 $[G]_s=[H]_s ~({\rm in}~\ccG)\Longleftrightarrow [C(G)]_s=[C(H)]_s~({\rm in}~\ccP) $.
\end{theorem}

\proof
1. The equivalence  is a direct consequence of 
Proposition \ref{prop_weak_and_s-dismant}.

2. \underline{$\Longrightarrow$} : Follows from 
Proposition \ref{prop_C_and_s-dismant}.

\underline{$\Longleftarrow$} : If $ [C(G)]_s=[C(H)]_s$, we get  $ [Comp(C(G))]_s=[Comp(C(H))]_s$
by the first assertion of the theorem.
As $Comp \circ C=Bd$, we have $[G']_s=[H']_s$ and the conclusion 
follows from Proposition  \ref{prop_subdivision}.
\endproof

\subsection{The triangle $(\mathscr{G},\ccP,\ccK)$}

Let us recall that in $\ccP$ there is also a notion of barycentric 
subdivision  $Bd :\ccP \to \ccP$ (for a poset $P$,
 $Bd(P)=P'$   is given by  the chains of  $P$ ordered by inclusion of underlying sets). 
There is also two classical applications,
$\Delta_{\ccP} :\ccP \to \ccK$ (the simplices of $\Delta_{\ccP}(P)$, 
the {\sl order complex} of $P$,
are   the $<x_0,x_1,\ldots,x_n>$  for every chain $x_0<x_1<\ldots <x_n$ of $P$) 
and $\Pi : \ccK \to \ccP$
(the elements of $\Pi(K)$, the {\sl face poset} of $K$,  are the simplices of  $K$
ordered by inclusion). Thus, we get the triangle $(\mathscr{G},\ccP,\ccK)$  
given in Fig. \ref{GPK}.

\begin{figure}[h]
 \begin{center}
\psset{unit=0.6 cm}
\pspicture(-2,-0.5)(6,6.5)
\rput(0,0){\large $\mathscr{G}$}
\rput(0,6){\large $\ccP$}
\rput(5,3){\large $\ccK$}
\psline{->}(1,5.1)(4,3.6)
\rput(2.5,5.5){\large $\Pi$}
\psline{->}(4,4.2)(1,5.7)
\rput(2.5,3.8){\large $\Delta_{\ccP}$}
\psline{->}(-0.3,4.8)(-0.3,1.2)
\rput(-1.4,3){\large $Comp$}
\psline{->}(0.3,1.2)(0.3,4.8)
\rput(1.1,3){\large $C$}
\psline{->}(1,0.9)(4,2.4)
\rput(2.5,2.3){\large $\Delta_{\mathscr{G}}$}
\psline{->}(4,1.8)(1,0.3)
\rput(2.5,0.5){\large $\Gamma$}
\endpspicture
\caption{The triangle $(\mathscr{G},\ccP,\ccK)$}
\label{GPK}
\end{center}
\end{figure}

\noindent Let us list some easy properties of  this 
triangle\footnote{In fact, we can consider $\ccG$, $\ccK$ and $\ccP$ 
as categories (with obvious morphisms) and it is easy to verify that all 
applications in the triangle $(\mathscr{G},\ccP,\ccK)$
are covariant functors. The reader not acquainted with the notions of functors 
or categories may refer to book \cite{kozlov_CAT}.}:
\ms

\begin{proposition} 
 1. $\Pi \circ \Delta_{\ccP}=C \circ Comp=Bd$ (in $\ccP$), 
$\Delta_{\ccP} \circ \Pi=\Delta_{\mathscr{G}}\circ \Gamma=Bd$
(in $\ccK$), $Comp \circ C=\Gamma \circ \Delta_{\mathscr{G}}=Bd$ (in $\mathscr{G}$).

2. We have the `` commutative triangles'' : 
$\Delta_{\ccP}=\Delta_{\mathscr{G}} \circ Comp$, 
$C=\Pi \circ \Delta_{\mathscr{G}} $ and
$\Gamma=Comp \circ \Pi$.

3. We have the `` commutative triangles  up to subdivision'' :
$\Gamma \circ \Delta_{\ccP}=Bd \circ Comp$, 
$\Delta_{\ccP} \circ C=Bd \circ \Delta_{\mathscr{G}}$ 
et $C \circ \Gamma=Bd \circ \Pi$. 
\end{proposition}

\noindent Now from Theorem \ref{theoremGK} and Theorem \ref{theoremGP}, 
we get another proof
of the Theorem (part of \cite[First main Theorem 3.9]{barmin}):

\begin{theorem}\label{theoremPK}
1. Let $P,Q \in \ccP$. Then $P$ and $Q$ are simply equivalent  if, and only if, 
$\Delta_{\ccP}(P)$ and $\Delta_{\ccP}(Q)$ have the same simple-homotopy type.

2. Let $K,L \in \ccK$. Then $K$ and $L$ have the same simple-homotopy type  
if, and only if, 
$\Pi(K)$ and $\Pi(L)$ are simply equivalent.
\end{theorem}

\begin{remark}
The image of $\Delta_{\ccG}$ is exactly the set of flag complexes 
but not all flag complexes are in 
 the image of $\Delta_{\ccP}$ (for example, the cyclic graph $C_5$ with 5 
vertices may be considered as a flag
complex and is not in the image of $\Delta_{\ccP}$; equivalently, $C_5$ 
is not a comparability graph).
\end{remark}

\section{The weak-s-dismantlability}\label{sec_weak}

\begin{definition}\label{s_dismantlable_edge}
Let $G \in \ccG$. An edge $gg'$ of $G$ will be called {\sl s-dismantlable} if
 $N_G(g)\cap N_G(g')$ is nonempty and dismantlable.
\end{definition}

We shall say that $G \redws H$ if we can go from $G$ to $H$ either by deleting s-dismantlable 
vertices or by deleting s-dismantlable edges. 

\begin{definition} \label{def_ws_homotopy_type}
Two graphs $G$ and $H$ have the same ws-homotopy type
if there is a sequence  $G=J_1,\ldots,J_k=H$ in  $\mathscr{G}$ such that
$G =J_1\stackrel{ws}{\to} J_2 \stackrel{ws}{\to}  \ldots \stackrel{ws}{\to} J_{k-1} 
\stackrel{ws}{\to} J_k=H$ 
where each arrow $\stackrel{ws}{\to} $ represents either 
the suppression or the addition of an s-dismantlable vertex, or 
the suppression or the addition of an s-dismantlable edge.
\end{definition}

\noindent  This defines an equivalence relation in $\mathscr{G}$
 and we shall denote by $[G]_{ws}$
the equivalence class representing the ws-homotopy type of a graph $G$. 
Of course,  $G \red H$ implies $G \redws H$ and the example 
given in remark \ref{s_converse_G_to_K_false}
shows that the reverse implication is false in general;
nevertheless, s-homotopy type and ws-homotopy type are the same:

\begin{proposition} For every $G \in \mathscr{G}$, we have $[G]_s=[G]_{ws}$.
\end{proposition}

\proof   The inclusion $[G]_s \subseteq [G]_{ws}$ follows 
from  $G \red H \Longrightarrow G \redws H$.
Now, we have seen (Lemma \ref{edge_s_deletable}) 
that the deletion of an s-dismantlable edge 
corresponds to the addition of an s-dismantlable vertex followed by the suppression of an
s-dismantlable vertex. This means that a sequence  
$G =J_1\stackrel{ws}{\to} J_2 \stackrel{ws}{\to}  \ldots \stackrel{ws}{\to} J_{k-1} 
\stackrel{ws}{\to} J_k=H$ can be rewritten as a sequence from $G$ to $H$ using 
only suppressions and additions of s-dismantlable vertices and proves that 
$[G]_s \supseteq [G]_{ws}$.
\endproof

Actually, the weak-s-dismantlability behaves well with the map 
$\Delta_{\mathscr{G}}$. The 1-skeleton of a simplicial complex 
can be considered as a graph
(whose vertices are given by the 0-simplices and the edges are given by the 1-simplices). 
Following the notation of \cite{cyy}, this defines a map  $sk :\mathscr{K} \to \mathscr{G}$
(if $K$ a simplicial complex, $sk(K)$ is its 1-skeleton taken as a graph). We note
that $sk(\Delta_{\mathscr{G}}(G))=G$ for all $G\in \mathscr{G}$. We have:

\begin{lemma}\label{lemma_sk}
Let $K$ be a flag simplicial complex and $L=K\setminus \{\sigma,\tau\}$.
If $ K\red L$ is an elementary simplicial collapse
(\textsl{i.e.}, if $\tau$ and $\s$ are faces of $K$ such that 
$\sigma$ is the unique face strictly containing $\tau$), then
$ sk(K) \redws sk(L)~{\rm or}~sk(K)=sk(L)$.
\end{lemma}

\proof The simplicial collapse $K \red L$ says that we obtain $L$ by deleting 
successively various  pairs of simplices $\{\sigma,\tau\}$ where $\tau$ is a free  
face of $\sigma$, so it is sufficient to prove that 
$sk(K) \redws sk(K\setminus \{\sigma,\tau\})$
for an elementary collapse $K \red (K\setminus \{\sigma,\tau\})$. 
Of course, if $\tau$ is a $k$-simplex, 
then $\sigma$ is a $(k+1)$-simplex and we consider the three cases
$k=0$, $k=1$ and $k\geq 2$. If $k=0$, $\tau=<a>$ is a vertex 
(or 0-simplex) belonging to a unique 1-simplex $\sigma=<a,b>$. 
In this case, we have $sk(K) \red (sk(K)\setminus a)$ because $a$ is a vertex 
dominated by the vertex $b$ in $sk(K)$ and 
$sk(K) \red sk(K\setminus \{\sigma,\tau\})=sk(K)\setminus a$ (and also 
$sk(K) \redws sk(K\setminus \{\sigma,\tau\})$. If $k=1$, $\tau=<a,b>$ is a 1-simplex 
such that there is a unique vertex $c$ such that  $<a,b,c>$ 
is a 2-simplex (named $\sigma$).
It follows that $c$ is the unique vertex of $sk(K)$ adjacent to $a$ and $b$; in other terms,
$N_{sk(K)}(a)\cap N_{sk(K)}(b)$ is reduced to the vertex $c$ and this shows that
the edge $ab$ is s-dismantlable in $sk(K)$. Now, it is clear that $sk(K\setminus
\{\sigma,\tau\})=sk(K\setminus \{<a,b,c>,<a,b>\})=sk(K\setminus \{<a,b>\})
=sk(K)\setminus ab$, so we obtain $sk(K) \redws sk(K\setminus \{\sigma,\tau\})
=sk(K)\setminus ab$. Finally, if $k\geq 2$, the suppression of the pair 
$\{\sigma,\tau\}$ in $K$ does not affect the 1-skeleton of $K$, \textit{i.e.} 
 $sk(K) = sk(K\setminus \{\sigma,\tau\})$.
\endproof

\begin{proposition}\label{ws_prop_G_to_K} Let $G,H \in \mathscr{G}$. 
Then,  $G \redws H \Longrightarrow \Delta_{\mathscr{G}}(G) \red \Delta_{\mathscr{G}}(H)$.
\end{proposition}

\proof By replacing  the 0-simplex  $<g>$ 
by the 1-simplex $<g,g'>$ in the  proof of Proposition \ref{s_prop_G_to_K},
we get  $\Delta_{\mathscr{G}}(G) \red \Delta_{\mathscr{G}}(G\setminus gg')$
when the edge $gg'$ is {\sl s-dismantlable} in $G$. This shows that 
$G \redws H$ implies $\Delta_{\mathscr{G}}(G) \red \Delta_{\mathscr{G}}(H)$.
\endproof

It is important to note 
that we can find a graph $G$ whose vertices 
and edges are all non s-dismantlable and such that
$\Delta_{\mathscr{G}}(G)$ collapses on a strict subcomplex 
which does not admit any collapsible pair and which is not a flag subcomplex; 
the 6-regular graph given in appendix 
provides such an example.

\section{Relation with graph homotopy of Chen, Yau and Yeh}\label{sec_Chen}

In \cite{iva}, Ivashchenko introduces the notion of {\sl contractible transformations} 
and calls {\sl contractible}   the trivial graph 
(the graph reduced to a point) and every graph
obtained from the trivial graph by application
of these contractible transformations.
In what follows, to avoid any confusion, 
we call {\sl I-contractibility}  the contractibility in the sense of Ivashchenko.
In a graph $G$, the contractible transformations are the deletion of a vertex $g$ if 
$N_G(g)$ is I-contractible, the deletion of an edge $gg'$ if $N_G(g)\cap N_G(g')$ 
is I-contractible,
 the addition of a vertex $x$ if $N_{G \cup x}(x)$ is I-contractible and the addition 
of an edge
between $g$ and $g'$ if $g \not\sim g'$ and $N_G(g)\cap N_G(g')$ is I-contractible.
 From these operations, in \cite{cyy}, the authors introduce the {\sl graph homotopy type}
of a graph $G$ that we shall call here {\sl $I$-homotopy type}. 
Let us  say that a vertex $g$ is I-dismantlable if  $N_G(g)$ is I-contractible.
The \cite[Lemma 3.4]{cyy} shows that we can reduce the four operations above
to the two operations of deletion or addition of  I-dismantlable vertices. 
Thus, one can define $G \redi H$ as the passage from $G$ to $H$ by suppression
of I-dismantlable vertices and $G \expi H$ as the passage from $G$ to $H$ by addition
of I-dismantlable vertices. From this, we get the $I$-equivalence class of a graph
(similarly to definition \ref{def_s_homotopy_type}, we say that
two graphs $G$ and $H$ have the same I-homotopy type
if there is a sequence  $G=J_1,\ldots,J_k=H$ in  $\mathscr{G}$ such that
$G =J_1\stackrel{I}{\to} J_2 \stackrel{I}{\to}  \ldots \stackrel{I}{\to} J_{k-1} 
\stackrel{I}{\to} J_k=H$ 
where each arrow $\stackrel{I}{\to} $ represents the suppression of a 
I-dismantlable vertex
or the addition of a I-dismantlable vertex). So, $[G]_I$ denotes the I-homotopy 
type of a graph $G$,
\textit{i.e.} the graph homotopy type of $G$ in the terminology of \cite{cyy}.  
In that way, $G$ is  I-contractible if, and only if,
 $[G]_I=[pt]_I$.
\ms

\begin{proposition} \label{lien_avec_cyy}
Let  $g\in G$.

1. If $g$ is s-dismantlable, then $g$ is $I$-dismantlable. 

2. If $G$ is s-dismantlable (\textit{i.e.} $[G]_s =[pt]_s$),  then $G$ is I-contractible.
\end{proposition}

\proof  1. Suppose that   $g$ is s-dismantlable in $G$;
then $N_G(g)$ is dismantlable.
But from the definition of I-contractibility, it is clear that  
``dismantlable $\Longrightarrow$ I-contractible'' 
(the open neighborhood of a dismantlable vertex is a cone and this 
proves that a dismantlable
vertex is I-dismantlable);
 thus $N_G(g)$ is I-contractible, i.e, $g$ is I-dismantlable. 

2. It is a consequence of the assertion 1.
\endproof

It also follows from the assertion 1 of  Proposition  \ref{lien_avec_cyy} 
that if two graphs $G$ and $H$
have the same s-homotopy type, then they have the same I-homotopy type. We 
are unaware if the converse is true:
\vspace{2 mm}

\noindent {\bf Question:} \emph{Let $G \in \mathscr{G}$.
Are the s-homotopy type of $G$ and the $I$-homotopy type of $G$  identical~?}

\appendix

\section*{Appendix: A particular non ws-reducible graph}

\begin{figure}
\psset{unit=0.2 cm}
\pspicture(-30,-18)(25,24)
\rput(0,16){$\bullet$}
\uput[u](0,16){\large $1$}
\rput(0,-16){$\bullet$}
\uput[d](0,-16){\large $2$}
\rput(6,0){$\bullet$}
\rput(4.5,0.8){ $6$}
\rput(9,0){$\bullet$}
\uput[r](9,0){\large $x$}
\rput(24,0){$\bullet$}
\uput[r](24,0){\large $3$}
\rput(11,3){$\bullet$}
\uput[ur](11,3){ $5$}
\rput(11,-3){$\bullet$}
\uput[dr](11,-3){ $4$}
\rput(-6,0){$\bullet$}
\rput(-4.5,0.8){ $9$}
\rput(-24,0){$\bullet$}
\uput[l](-24,0){ $3$}
\rput(-11,3){$\bullet$}
\uput[ul](-11,3){ $8$}
\rput(-11,-3){$\bullet$}
\uput[dl](-11,-3){ $7$}
\psline(0,16)(6,0)(11,3)(0,16)(24,0)(11,3)(11,-3)(24,0)(0,-16)(11,-3)(6,0)(0,-16)(0,16)     
\psline(0,16)(-6,0)(-11,3)(0,16)(-24,0)(-11,3)(-11,-3)(-24,0)(0,-16)(-11,-3)(-6,0)(0,-16)
\psline(6,0)(9,0)
\psline(11,3)(9,0)
\psline(11,-3)(9,0)
\psline(-6,0)(6,0)
\psline(-11,3)(11,3)
\psline(-11,-3)(11,-3)
\pscurve(9,0)(0,-2)(-6,0)
\pscurve(9,0)(0,-5)(-11,-3)
\pscurve(9,0)(0,5)(-11,3)
\rput(-20,15){Figure 8.a.: }
\rput(-20,12){Graph $G$}
\rput(36,2){$\bullet$}
\rput(46,2){$\bullet$}
\rput(36,12){$\bullet$}
\rput(46,12){$\bullet$}
\rput(41,5){$\bullet$}
\rput(41,9){$\bullet$}
\psline(36,2)(46,2)(46,12)(36,12)(36,2)
\psline(36,2)(41,5)(46,2)
\psline(36,12)(41,9)(46,12)
\psline(41,9)(41,5)
\rput(41,-2){Figure 8.b.}
\rput(36,-10){$\bullet$}
\rput(46,-10){$\bullet$}
\rput(41,-10){$\bullet$}
\psline(36,-10)(41,-10)
\rput(41,-14){Figure 8.c}
\endpspicture
\caption{$\Delta_{\ccG}(G)$ is non collapsible on a flag subcomplex}
\label{non ws-reducible}
\end{figure}

\noindent Some properties of the graph $G$ of Fig. \ref{non ws-reducible}: \vspace{ 2 mm}

\begin{enumerate}
\item $G$ is a 6-regular graph with $\vert V(G)\vert =10$ and $\vert E(G)\vert =30$.
\item For every vertex $a$ of $G$, the subgraph $N_G(a)$ is not dismantlable because
it is isomorphic to the graph shown in figure 8.b. 
\item For every edge $ab$ of $G$, the subgraph $N_G(a)\cap N_G(b)$ is  
not dismantlable because
it is isomorphic to the (disconnected) graph shown in figure 8.c.  
So there is no s-dismantlable edge in $G$.
\item The cliques (or {\sl maximal complete subgraphs}) of $G$ are of order 4 or 3. 
The cliques of order 3 are 
$$[1,2,3]~~~[1,5,6]~~~[2,4,6]~~~[1,8,9]~~~[2,7,9]~~~[3,4,5]~~~[3,7,8]~~~[x,5,8]~~~[x,6,9]
~~~~[x,4,7]$$
and the cliques of order 4  are:
$$c_1=[1,3,5,8]~~~~~c_2=[1,2,6,9]~~~~~c_3=[2,3,4,7]~~~~~c_4=[x,4,5,6]~~~~c_5=[x,7,8,9]$$
\item In the simplicial complex $\Delta_{\ccG}(G)$, there are five tetrahedras $\sigma_i$
 corresponding to the five 4-cliques $c_i$.
Corresponding to each tetrahedron $\sigma_i$, 
each pair  $(\sigma_i,\tau)$ (with $\tau$ being any 
maximal proper face of $\sigma_i$) is a collapsible pair.
\item Let $K$ be a subcomplex obtained from $\Delta_{\ccG}(G)$ after collapsing the five 
tetrahedras $\tau_i$, $1\leq i\leq 5$. 
\begin{itemize}
\item There is no collapsible pair in $K$. Indeed, a collapsible pair in $K$
must be of the form $(\sigma',\tau')$ with $\sigma'$ a triangle (or 2-simplex) 
and $\tau'$ an edge  of 
$\sigma'$ which is not the edge of another triangle. 
But we see from the lists or 3-cliques and 4-cliques that every edge
of $G$  appears exactly once in the list of 3-cliques and
 exactly once in the list of 4-cliques; so, even after removing 5 
collapsible pairs corresponding to the 5 tetrahedras,
every edge appears in at least two triangles (and is not a free edge).
\item $K$ is not a flag complex. For example, let $(\sigma,\tau)$ be a pair 
which has been collapsed in 
$\Delta_{\ccG}(G)$: $\sigma$ is a tetrahedron and $\tau$, a maximal 
proper face of $\sigma$, is 
of the form $<a,b,c>$. So $\tau$ is a non-simplex of $K$ with 3 vertices 
and every face of $\tau$ is a simplex of $K$.
\end{itemize}
\end{enumerate}

\noindent In conclusion: 

\begin{itemize}
\item $G$ is a non ws-reducible graph (\textit{i.e.}, 
there is not a strict subgraph $H$ of $G$
such that $G \redws H$).
\item We can find a strict subcomplex $K$ of $\Delta_{\ccG}(G)$ 
such that $\Delta_{\ccG}(G) \red K$.
\item Every strict subcomplex $K$ of $\Delta_{\ccG}(G)$ such 
that $\Delta_{\ccG}(G) \red K$ is not a flag complex.
\end{itemize}

\vspace{5 mm}

\noindent Romain Boulet, Universit$\acute{\rm e}$
Toulouse Le Mirail, 5, All$\acute{\rm e}$es Antonio Machado, 
31058 Toulouse Cedex 9;
email:boulet@univ-tlse2.fr
\vspace{2 mm}

\noindent Etienne Fieux, Universit$\acute{\rm e}$ Paul Sabatier, 
118, Route de Narbonne - F-31062 Toulouse cedex 9 ; 
email: fieux@math.univ-toulouse.fr
\vspace{2 mm}

\noindent Bertrand Jouve, Universit$\acute{\rm e}$ Toulouse Le Mirail, 5, All$\acute{\rm e}$es Antonio Machado, 
31058 Toulouse Cedex 9;
email: jouve@univ-tlse2.fr

\end{document}